\documentclass[11pt]{article}
\usepackage{a4wide}

\newcommand{\ba}{\begin{array}}
\newcommand{\ea}{\end{array}}
\newcommand{\be}{\begin{equation}}
\newcommand{\ee}{\end{equation}}
\newcommand{\la}{\label}
\newcommand{\bea}{\begin{eqnarray}}
\newcommand{\eea}{\end{eqnarray}}
\newcommand{\ch}{\choose}
\renewcommand{\l}{\left}
\renewcommand{\r}{\right}
\newcommand{\n}{\nonumber}
\newcommand{\nn}{\nonumber \\}
\newcommand{\ds}{\displaystyle}
\newcommand{\ndots}{n=0,1,2,\ldots}
\renewcommand{\a}{\alpha}
\renewcommand{\b}{\beta}
\newcommand{\G}{\Gamma}
\renewcommand{\L}{L_n^{(\a)}(x)}
\newcommand{\KL}{L_n^{\a,M}(x)}
\newcommand{\GL}{L_n^{\a,M,N}(x)}
\renewcommand{\P}{P_n^{(\a,\b)}(x)}
\newcommand{\GP}{P_n^{\a,\b,M,N}(x)}
\newcommand{\SP}{P_n^{(\a,\a)}(x)}
\newcommand{\SGP}{P_n^{\a,\a,M,M}(x)}
\newcommand{\set}[1]{\left\{#1\right\}_{n=0}^{\infty}}
\newcommand{\hyp}[5]{\mbox{}_{#1}F_{#2}
\left(\left.{#3 \atop #4}\right|#5\right)}

\begin{document}

\title{The Jacobi inversion formula}
\author{J.~Koekoek \and R.~Koekoek}
\date{ }
\maketitle

\begin{abstract}
We look for differential equations satisfied by the generalized Jacobi
polynomials $\set{\GP}$ which are orthogonal on the interval $[-1,1]$ with
respect to the weight function
$$\frac{\G(\a+\b+2)}{2^{\a+\b+1}\G(\a+1)\G(\b+1)}(1-x)^{\a}(1+x)^{\b}+
M\delta(x+1)+N\delta(x-1),$$
where $\a>-1$, $\b>-1$, $M\ge 0$ and $N\ge 0$.

In order to find explicit formulas for the coefficients of these
differential equations we have to solve systems of equations of the form
$$\sum_{i=1}^{\infty}A_i(x)D^i\P=F_n(x),\;n=1,2,3,\ldots,$$
where the coefficients $\l\{A_i(x)\r\}_{i=1}^{\infty}$ are independent of
$n$. This system of equations has a unique solution given by
$$A_i(x)=2^i\sum_{j=1}^i\frac{\a+\b+2j+1}{(\a+\b+j+1)_{i+1}}
P_{i-j}^{(-\a-i-1,-\b-i-1)}(x)F_j(x),\;i=1,2,3,\ldots.$$
This is a consequence of the Jacobi inversion formula
\bea & &\sum_{k=j}^i\frac{\a+\b+2k+1}{(\a+\b+k+j+1)_{i-j+1}}\times{}\nn
& &{}\hspace{1cm}{}\times
P_{i-k}^{(-\a-i-1,-\b-i-1)}(x)P_{k-j}^{(\a+j,\b+j)}(x)=\delta_{ij},
\;j\le i,\;i,j=0,1,2,\ldots,\n\eea
which is proved in this paper.
\end{abstract}

\vfill

\begin{tabular}{ll}
Keywords : & Differential equations, Generalized Jacobi polynomials,\\
& Inversion formulas.\\[5mm]
\multicolumn{2}{l}{1991 Mathematics Subject Classification :
Primary 33C45 ; Secondary 34A35.}
\end{tabular}

\newpage

\section{Introduction}

In \cite{Koorn} T.H.~Koornwinder introduced the generalized Jacobi
polynomials $\set{\GP}$ which are orthogonal on the interval $[-1,1]$ with
respect to the weight function
$$\frac{\G(\a+\b+2)}{2^{\a+\b+1}\G(\a+1)\G(\b+1)}(1-x)^{\a}(1+x)^{\b}+
M\delta(x+1)+N\delta(x-1),$$
where $\a>-1$, $\b>-1$, $M\ge 0$ and $N\ge 0$. As a limit case he also found
the generalized Laguerre polynomials $\set{\KL}$ which are orthogonal on the
interval $[0,\infty)$ with respect to the weight function
$$\frac{1}{\G(\a+1)}x^{\a}e^{-x}+M\delta(x),$$
where $\a>-1$ and $M\ge 0$. These generalized Jacobi polynomials and
generalized Laguerre polynomials are related by the limit
$$\KL=\lim_{\b\rightarrow\infty}P_n^{\a,\b,0,M}\l(1-\frac{2x}{\b}\r).$$

In \cite{DV} we proved that for $M>0$ the generalized Laguerre polynomials
satisfy a unique differential equation of the form
$$M\sum_{i=0}^{\infty}a_i(x)y^{(i)}(x)+xy''(x)+(\a+1-x)y'(x)+ny(x)=0,$$
where $\l\{a_i(x)\r\}_{i=0}^{\infty}$ are continuous functions on the real
line and $\l\{a_i(x)\r\}_{i=1}^{\infty}$ are independent of the degree $n$.
In \cite{Bav} H.~Bavinck found a new method to obtain the main result of
\cite{DV}. This inversion method was found in a similar way as was done in
\cite{Char} in the case of generalizations of the Charlier polynomials. See
also section~4 for more details. In \cite{Soblag} we used this inversion
method to find all differential equations of the form
\bea\la{DVLag}& &M\sum_{i=0}^{\infty}a_i(x)y^{(i)}(x)+
N\sum_{i=0}^{\infty}b_i(x)y^{(i)}(x)+{}\nn
& &\hspace{1cm}{}+MN\sum_{i=0}^{\infty}c_i(x)y^{(i)}(x)+
xy''(x)+(\a+1-x)y'(x)+ny(x)=0,\eea
where the coefficients $\l\{a_i(x)\r\}_{i=1}^{\infty}$,
$\l\{b_i(x)\r\}_{i=1}^{\infty}$ and $\l\{c_i(x)\r\}_{i=1}^{\infty}$ are
independent of $n$ and the coefficients $a_0(x)$, $b_0(x)$ and $c_0(x)$ are
independent of $x$, satisfied by the Sobolev-type Laguerre polynomials
$\set{\GL}$ which are orthogonal with respect to the inner product
$$<f,g>\;=\frac{1}{\G(\a+1)}\int_0^{\infty}x^{\a}e^{-x}f(x)g(x)dx+Mf(0)g(0)+
Nf'(0)g'(0),$$
where $\a>-1$, $M\ge 0$ and $N\ge 0$. These Sobolev-type Laguerre
polynomials $\set{\GL}$ are generalizations of the generalized Laguerre
polynomials $\set{\KL}$. In fact we have
$$L_n^{\a,M,0}(x)=\KL\;\mbox{ and }\;L_n^{\a,0}(x)=\L.$$

In this paper we will prove an inversion formula involving the classical
Jacobi polynomials which can be used to find differential equations of the
form
\bea\la{DVJac}& &M\sum_{i=0}^{\infty}a_i(x)y^{(i)}(x)+
N\sum_{i=0}^{\infty}b_i(x)y^{(i)}(x)+
MN\sum_{i=0}^{\infty}c_i(x)y^{(i)}(x)+{}\nn
& &{}\hspace{1cm}
{}+(1-x^2)y''(x)+\l[\b-\a-(\a+\b+2)x\r]y'(x)+n(n+\a+\b+1)y(x)=0,\eea
where the coefficients $\l\{a_i(x)\r\}_{i=1}^{\infty}$,
$\l\{b_i(x)\r\}_{i=1}^{\infty}$ and $\l\{c_i(x)\r\}_{i=1}^{\infty}$ are
independent of $n$ and the coefficients $a_0(x)$, $b_0(x)$ and $c_0(x)$ are
independent of $x$, satisfied by the generalized Jacobi polynomials
$\set{\GP}$. In \cite{Madrid} we applied the special case $\b=\a$ of this
inversion formula to solve the systems of equations obtained in
\cite{Symjac}.

The inversion formula for the Charlier polynomials obtained in \cite{Char}
(see also section~4) was also used in \cite{SobChar} to find difference
operators with Sobolev-type Charlier polynomials as eigenfunctions.

In \cite{Meixner} H.~Bavinck and H.~van~Haeringen used similar inversion
formulas to find difference equations for generalized Meixner polynomials.

\section{The classical Laguerre and Jacobi polynomials}

In this section we list the definitions and some properties of the classical
Laguerre and Jacobi polynomials which we will use in this paper. For details
the reader is referred to \cite{Chihara}, \cite{AS} and \cite{Szego}.

The classical Laguerre polynomials $\set{\L}$ can be defined by
\be\la{defLag}\L=\sum_{k=0}^n(-1)^k{n+\a \ch n-k}\frac{x^k}{k!},\;\ndots\ee
for all $\a$. Their generating function is given by
\be\la{genLag}(1-t)^{-\a-1}\exp\l(\frac{xt}{t-1}\r)=
\sum_{n=0}^{\infty}\L t^n\ee
and for all $n\in\{0,1,2,\ldots\}$ we have
\be\la{diffLag}D^i\L=(-1)^iL_{n-i}^{(\a+i)}(x),\;i=0,1,2,\ldots,n,\ee
where $D=\ds\frac{d}{dx}$ denotes the differentiation operator. The Laguerre
polynomials satisfy the linear second order differential equation
\be\la{dvLag}xy''(x)+(\a+1-x)y'(x)+ny(x)=0.\ee
It is well-known that
\be\la{formLag}\frac{x^n}{n!}=
\sum_{k=0}^n(-1)^k{n+\a \ch n-k}L_k^{(\a)}(x),\;\ndots.\ee
This formula can easily be proved by using definition (\ref{defLag}) and
changing the order of summation as follows
\bea\sum_{k=0}^n(-1)^k{n+\a \ch n-k}L_k^{(\a)}(x)&=&\sum_{k=0}^n(-1)^k
{n+\a \ch n-k}\sum_{j=0}^k(-1)^j{k+\a \ch k-j}\frac{x^j}{j!}\nn
&=&\sum_{j=0}^n\sum_{k=0}^{n-j}(-1)^k{n+\a \ch n-j-k}{j+k+\a \ch k}
\frac{x^j}{j!}\nn
&=&\sum_{j=0}^n{n+\a \ch n-j}\frac{x^j}{j!}\sum_{k=0}^{n-j}(-1)^k{n-j \ch k}
=\frac{x^n}{n!},\;\ndots.\n\eea

The classical Jacobi polynomials $\set{\P}$ can be defined by
\bea\la{defJac1}\P&=&\sum_{k=0}^n\frac{(n+\a+\b+1)_k}{k!}
\frac{(\a+k+1)_{n-k}}{(n-k)!}\l(\frac{x-1}{2}\r)^k,\;\ndots\\
&=&\la{defJac2}(-1)^n\sum_{k=0}^n\frac{(-n-k-\a-\b)_k}{k!}
\frac{(-n-\a)_{n-k}}{(n-k)!}\l(\frac{x-1}{2}\r)^k,\;\ndots\\
&=&\la{defJac3}2^{-n}\sum_{k=0}^n{n+\a \ch n-k}{n+\b \ch k}
(x-1)^k(x+1)^{n-k},\;\ndots\eea
for all $\a$ and $\b$. For all $n\in\{0,1,2,\ldots\}$ we have
\be\la{diffJac}D^i\P=\frac{(n+\a+\b+1)_i}{2^i}P_{n-i}^{(\a+i,\b+i)}(x),
\;i=0,1,2,\ldots,n.\ee
These Jacobi polynomials satisfy the linear second
order differential equation
\be\la{dvJac}(1-x^2)y''(x)+\l[\b-\a-(\a+\b+2)x\r]y'(x)
+n(n+\a+\b+1)y(x)=0.\ee
Further we have for $\a+\b+1>0$ (compare with \cite{Luke}, page 277, formula
(30))
\be\la{formJac*}\l(\frac{1-x}{2}\r)^n=\sum_{k=0}^n
\frac{(-n)_k(\a+k+1)_{n-k}(\a+\b+2k+1)}{(\a+\b+k+1)_{n+1}}P_k^{(\a,\b)}(x),
\;\ndots.\ee
This formula is much less known than formula (\ref{formLag}) for the
Laguerre polynomials. However, the proof is quite similar. In section~5 we
will prove a much more general formula.

We remark that (\ref{formJac*}) can be written in a more general form as
\bea\la{formJac} & &\sum_{k=0}^n
\frac{(-n)_k(\a+\b+1)_k(\a+k+1)_{n-k}(\a+\b+2k+1)}{\G(\a+\b+n+k+2)}
P_k^{(\a,\b)}(x)\nn
&=&\frac{1}{\G(\a+\b+1)}\l(\frac{1-x}{2}\r)^n,\;\ndots,\eea
which is valid for all $\a$ and $\b$.

The Jacobi polynomials $\set{\P}$ and the Laguerre polynomials $\set{\L}$
are related by the limit
\be\la{limit}\L=\lim_{\b\rightarrow\infty}
P_n^{(\a,\b)}\l(1-\frac{2x}{\b}\r).\ee
We remark that if we replace $x$ by $\ds 1-\frac{2x}{\b}$ in (\ref{formJac*}),
multiply by $\b^n$ and let $\b$ tend to infinity in the complex plane along
the halfline where $\a+\b$ is real and $\a+\b+1>0$ we obtain (\ref{formLag})
by using (\ref{limit}) and the fact that we have for all
$n\in\{0,1,2,\ldots\}$
$$(-n)_k(\a+k+1)_{n-k}=(-1)^k{n+\a \ch n-k}n!,k=0,1,2,\ldots,n.$$

\section{The systems of equations}

Let $\a>-1$. The Sobolev-type Laguerre polynomials $\set{\GL}$ can be
written as
$$\GL=A_0\L+A_1D\L+A_2D^2\L,\;\ndots,$$
where the coefficients $A_0$, $A_1$ and $A_2$ are given by
$$\l\{\ba{l}\ds A_0=1+M{n+\a \ch n-1}+\frac{n(\a+2)-(\a+1)}{(\a+1)(\a+3)}
N{n+\a \ch n-2}+{}\\
\ds\hspace{5cm}{}+\frac{MN}{(\a+1)(\a+2)}{n+\a \ch n-1}{n+\a+1 \ch n-2}\\  \\
\ds A_1=M{n+\a \ch n}+\frac{n-1}{\a+1}N{n+\a \ch n-1}+
\frac{2MN}{(\a+1)^2}{n+\a \ch n}{n+\a+1 \ch n-2}\\  \\
\ds A_2=\frac{N}{\a+1}{n+\a \ch n-1}+\frac{MN}{(\a+1)^2}{n+\a \ch n}
{n+\a+1 \ch n-1}.\ea\r.$$
For details concerning these Sobolev-type Laguerre polynomials and their
definition the reader is referred to \cite{Thesis} and \cite{SIAM}. Since
the classical Laguerre polynomials $\set{\L}$ satisfy the differential
equation (\ref{dvLag}) it is quite reasonable to look for differential
equations of the form (\ref{DVLag}) for these Sobolev-type Laguerre
polynomials $\set{\GL}$ in view of this definition and the fact that
$L_n^{\a,0,0}(x)=\L$. In \cite{Soblag} it is shown that this leads to eight
systems of equations for the coefficients $\l\{a_i(x)\r\}_{i=0}^{\infty}$,
$\l\{b_i(x)\r\}_{i=0}^{\infty}$ and $\l\{c_i(x)\r\}_{i=0}^{\infty}$. In
order to find these coefficients we have to solve systems of equations which
are of the form
$$\sum_{i=1}^{\infty}A_i(x)D^{i+k}\L=F_n(x),\;n=k+1,k+2,k+3,\ldots,$$
where $k\in\{0,1,2,\ldots\}$ and the coefficients
$\l\{A_i(x)\r\}_{i=1}^{\infty}$ are independent of $n$. In \cite{Soblag} it
is pointed out that this system of equations has a unique solution given by
$$A_i(x)=(-1)^{i+k}\sum_{j=1}^iL_{i-j}^{(-\a-i-k-1)}(-x)F_{j+k}(x),
\;i=1,2,3,\ldots.$$
This is an easy consequence of the Laguerre inversion formula
\be\la{invLag}\sum_{k=j}^iL_{i-k}^{(-\a-i-1)}(-x)L_{k-j}^{(\a+j)}(x)=
\delta_{ij},\;j\le i,\;i,j=0,1,2,\ldots,\ee
which was found by H.~Bavinck in \cite{Bav}. For more details the reader is
referred to \cite{Bav} and \cite{Soblag}. See also section~4 of this paper.

Now we take $\a>-1$ and $\b>-1$. The generalized Jacobi polynomials
$\set{\GP}$ can be written as
$$\GP=A_0\P+\l[A_1(1-x)-A_2(1+x)\r]D\P,\;\ndots,$$
where the coefficients $A_0$, $A_1$ and $A_2$ are given by
$$\l\{\ba{l}\ds A_0=1+
M\frac{\ds {n+\b \ch n-1}{n+\a+\b+1 \ch n}}{\ds {n+\a \ch n}}+
N\frac{\ds {n+\a \ch n-1}{n+\a+\b+1 \ch n}}{\ds {n+\b \ch n}}+{}\\  \\
\ds\hspace{7cm}{}+
MN\frac{(\a+\b+2)^2}{(\a+1)(\b+1)}{n+\a+\b+1 \ch n-1}^2\\ \\
\ds A_1=\frac{M}{\a+\b+1}
\frac{\ds {n+\b \ch n}{n+\a+\b \ch n}}{\ds {n+\a \ch n}}
+\frac{MN}{\a+1}{n+\a+\b \ch n-1}{n+\a+\b+1 \ch n}\\ \\
\ds A_2=\frac{N}{\a+\b+1}
\frac{\ds {n+\a \ch n}{n+\a+\b \ch n}}{\ds {n+\b \ch n}}
+\frac{MN}{\b+1}{n+\a+\b \ch n-1}{n+\a+\b+1 \ch n}.\ea\r.$$
Here we used the same definition as in \cite{Koorn}, but in a slightly
different notation. The case $\a+\b+1=0$ must be understood by continuity.
In view of this definition and the fact that the classical Jacobi
polynomials $\set{\P}$ satisfy the differential equation (\ref{dvJac}) it
is quite natural to look for differential equations of the form
(\ref{DVJac}) satisfied by these generalized Jacobi polynomials as was
already pointed out in \cite{JCAM}. Again this leads to eight systems of
equations for the coefficients $\l\{a_i(x)\r\}_{i=0}^{\infty}$,
$\l\{b_i(x)\r\}_{i=0}^{\infty}$ and $\l\{c_i(x)\r\}_{i=0}^{\infty}$. In
order to find these coefficients we have to solve systems of equations which
are of the form
$$\sum_{i=1}^{\infty}A_i(x)D^i\P=F_n(x),\;n=1,2,3,\ldots,$$
where the coefficients $\l\{A_i(x)\r\}_{i=1}^{\infty}$ are independent of
$n$. This system of equations has a unique solution given by
$$A_i(x)=2^i\sum_{j=1}^i\frac{\a+\b+2j+1}{(\a+\b+j+1)_{i+1}}
P_{i-j}^{(-\a-i-1,-\b-i-1)}(x)F_j(x),\;i=1,2,3,\ldots.$$
This is a consequence of the Jacobi inversion formula
\bea\la{invJac}& &\sum_{k=j}^i
\frac{\a+\b+2k+1}{(\a+\b+k+j+1)_{i-j+1}}\times{}\nn
& &{}\hspace{1cm}{}\times
P_{i-k}^{(-\a-i-1,-\b-i-1)}(x)P_{k-j}^{(\a+j,\b+j)}(x)=\delta_{ij},
\;j\le i,\;i,j=0,1,2,\ldots,\eea
which will be proved in this paper. Again, the case $\a+\b+1=0$ must be
understood by continuity. We remark that if we replace $x$ by
$\ds 1-\frac{2x}{\b}$ in (\ref{invJac}), multiply by $\b^{i-j}$ and let $\b$
tend to infinity along the positive real axis we obtain the Laguerre
inversion formula (\ref{invLag}) by using (\ref{limit}).

In \cite{Symjac} we found all differential equations of the form
\be\la{DVSGP}M\sum_{i=0}^{\infty}a_i(x)y^{(i)}(x)+
(1-x^2)y''(x)-2(\a+1)xy'(x)+n(n+2\a+1)y(x)=0,\ee
where $\l\{a_i(x)\r\}_{i=0}^{\infty}$ are continuous functions on the real
line and $\l\{a_i(x)\r\}_{i=1}^{\infty}$ are independent of $n$, satisfied
by the symmetric generalized ultraspherical polynomials $\set{\SGP}$ defined
by
$$\SGP=C_0\SP-C_1xD\SP,\;\ndots,$$
where
$$\l\{\ba{l}\ds C_0=1+\frac{2Mn}{\a+1}{n+2\a+1 \ch n}
+4M^2{n+2\a+1 \ch n-1}^2\\
\\
\ds C_1=\frac{2M}{2\a+1}{n+2\a \ch n}
+\frac{2M^2}{\a+1}{n+2\a \ch n-1}{n+2\a+1 \ch n}.\ea\r.$$
We remark that these polynomials form a special case ($\b=\a$ and $N=M$) of
the generalized Jacobi polynomials $\set{\GP}$, but the differential
equation (\ref{DVSGP}) has a very special form without a $M^2$-part. This is
explained by the fact that
$$\l\{\ba{l}\ds C_0=\l[1+2M{n+2\a+1 \ch n-1}\r]^2\\
\\
\ds C_1=\frac{2M}{2\a+1}{n+2\a \ch n}\l[1+2M{n+2\a+1 \ch n-1}\r].\ea\r.$$
This implies that the generalized ultraspherical polynomials satisfy the
same differential equation as the polynomials $\set{Q_n^{\a,\a,M,M}(x)}$
defined by
\bea & &Q_n^{\a,\a,M,M}(x)=\l[1+2M{n+2\a+1 \ch n-1}\r]\SP+{}\nn
& &{}\hspace{5cm}{}-\frac{2M}{2\a+1}{n+2\a \ch n}xD\SP,\;\ndots.\n\eea
However, this differential equation will appear not to be a special case of
the differential equation of the form (\ref{DVJac}) for the generalized Jacobi
polynomials, since the $MN$-part will not vanish if we take $\b=\a$ and
$N=M$. We aim to give a proof of this in a future publication. In
\cite{Madrid} we applied the special case $\b=\a$ of the Jacobi inversion
formula (\ref{invJac}) to solve the systems of equations obtained in
\cite{Symjac}.

\section{The inversion formulas}

In \cite{Char} H.~Bavinck and R.~Koekoek found the following inversion
formula involving Charlier polynomials
\be\la{invChar}\sum_{k=j}^iC_{i-k}^{(-a)}(-x)C_{k-j}^{(a)}(x)=
\delta_{ij},\;j\le i,\;i,j=0,1,2,\ldots.\ee
This formula is an easy consequence of the generating function (see for
instance \cite{AS})
$$e^{-at}(1+t)^x=\sum_{n=0}^{\infty}C_n^{(a)}(x)t^n.$$
In fact we have
\bea 1&=&e^{-at}(1+t)^xe^{at}(1+t)^{-x}=
\sum_{k=0}^{\infty}C_k^{(a)}(x)t^k\sum_{m=0}^{\infty}C_m^{(-a)}(-x)t^m\nn
&=&\sum_{n=0}^{\infty}\l(\sum_{k=0}^nC_k^{(a)}(x)C_{n-k}^{(-a)}(-x)\r)t^n.
\n\eea
Hence
$$\sum_{k=0}^nC_k^{(a)}(x)C_{n-k}^{(-a)}(-x)=\l\{\ba{ll}1, & n=0\\ \\
0, & n=1,2,3,\ldots.\ea\r.$$
Now (\ref{invChar}) easily follows by taking $n=i-j$ and shifting the
summation index. This formula was also used in \cite{SobChar} to find
difference operators with Sobolev-type Charlier polynomials as
eigenfunctions. In \cite{Meixner} a similar formula involving Meixner
polynomials was used to find difference equations for generalized Meixner
polynomials.

Formula (\ref{invChar}) can be interpreted as follows. If we define the
matrix $T=(t_{ij})_{i,j=0}^n$ with entries
$$t_{ij}=\l\{\ba{ll}C_{i-j}^{(a)}(x), & j\le i\\ \\0, & j>i,\ea\r.$$
then this matrix $T$ is a triangular matrix with determinant $1$ and the
inverse $U$ of this matrix is given by $T^{-1}=U=(u_{ij})_{i,j=0}^n$ with
entries
$$u_{ij}=\l\{\ba{ll}C_{i-j}^{(-a)}(-x), & j\le i\\ \\0, & j>i.\ea\r.$$
Therefore we call (\ref{invChar}) an inversion formula.

In the same way we find by using the generating function (\ref{genLag})
for the Laguerre polynomials
\be\la{invLag*}\sum_{k=j}^iL_{i-k}^{(\a)}(x)L_{k-j}^{(-\a-2)}(-x)=
\delta_{ij},\;j\le i,\;i,j=0,1,2,\ldots.\ee
However, this formula cannot be used to solve systems of equations of the
form
$$\sum_{i=1}^{\infty}A_i(x)D^i\L=F_n(x),\;n=1,2,3,\ldots$$
in view of the parametershift in (\ref{diffLag}).

In \cite{Bav} H.~Bavinck used a slightly different method to find the
Laguerre inversion formula (\ref{invLag}) from the generating function
(\ref{genLag}) for the Laguerre polynomials. In fact we have
\bea (1-t)^{i-j-1}&=&(1-t)^{-\a-j-1}\exp\l(\frac{xt}{t-1}\r)
(1-t)^{\a+i}\exp\l(\frac{-xt}{t-1}\r)\nn
&=&\sum_{k=0}^{\infty}L_k^{(\a+j)}(x)t^k
\sum_{m=0}^{\infty}L_m^{(-\a-i-1)}(-x)t^m\nn
&=&\sum_{n=0}^{\infty}
\l(\sum_{k=0}^nL_k^{(\a+j)}(x)L_{n-k}^{(-\a-i-1)}(-x)\r)t^n.\n\eea
This implies, by comparing the coefficients of $t^{i-j}$ on both sides, that
$$\sum_{k=0}^{i-j}L_k^{(\a+j)}(x)L_{i-j-k}^{(-\a-i-1)}(-x)=\delta_{ij},
\;j\le i,\;i,j=0,1,2,\ldots,$$
which is equivalent to (\ref{invLag}).

Formula (\ref{invLag}) can be interpreted as follows. If we define the
matrix $T=(t_{ij})_{i,j=0}^n$ with entries
$$t_{ij}=\l\{\ba{ll}L_{i-j}^{(\a+j)}(x), & j\le i\\ \\0, & j>i,\ea\r.$$
then this matrix $T$ is a triangular matrix with determinant $1$ and the
inverse $U$ of this matrix is given by $T^{-1}=U=(u_{ij})_{i,j=0}^n$ with
entries
$$u_{ij}=\l\{\ba{ll}L_{i-j}^{(-\a-i-1)}(-x), & j\le i\\ \\0, & j>i.\ea\r.$$

In case of the Jacobi polynomials the above methods seem not to be
applicable. In that case we have to find the inverse of the matrix
$T=(t_{ij})_{i,j=0}^n$ with entries
$$t_{ij}=D^jP_i^{(\a,\b)}(x),\;i,j=0,1,2,\ldots,n.$$
This matrix $T$ is also triangular and by using (\ref{diffJac}) the diagonal
entries equal
$$t_{ii}=D^iP_i^{(\a,\b)}(x)=\frac{(i+\a+\b+1)_i}{2^i},\;i=0,1,2,\ldots,n.$$
This implies that the determinant of $T$ is nonzero for each $n$ iff
$-(\a+\b+2)\notin\{0,1,2,\ldots\}$. In that case $T$ is invertible and if
the inverse $U$ is given by $T^{-1}=U=(u_{ij})_{i,j=0}^n$ then we must have
$$u_{ii}=\frac{1}{t_{ii}}=\frac{2^i}{(\a+\b+i+1)_i},\;i=0,1,2,\ldots,n.$$
In the next section we will give a proof of the Jacobi inversion formula
(\ref{invJac}), which is equivalent to
$$u_{ij}=\l\{\ba{ll}\ds\frac{(\a+\b+2j+1)2^i}{(\a+\b+j+1)_{i+1}}
P_{i-j}^{(-\a-i-1,-\b-i-1)}(x), & j\le i\\ \\0, & j>i.\ea\r.$$

\section{Proof of the Jacobi inversion formula}

In this section we will prove that
\bea\la{algJac}& &\sum_{k=0}^n
\frac{(\a+\b+2k+1)(\a+\b+1)_k}{\G(\a+\b+n+k+2)}P_k^{(\a,\b)}(x)
P_{n-k}^{(-n-\a-1,-n-\b-1)}(y)\nn
&=&\frac{1}{\G(\a+\b+1)}\frac{1}{n!}\l(\frac{x-y}{2}\r)^n,
\;n=0,1,2,\ldots,\eea
which holds for all $\a$ and $\b$.

Note that (\ref{formJac}) is a special case of (\ref{algJac}) since
\bea P_{n-k}^{(-n-\a-1,-n-\b-1)}(1)&=&\frac{(-n-\a)_{n-k}}{(n-k)!}=
(-1)^{n-k}\frac{(\a+k+1)_{n-k}}{(n-k)!}\nn
&=&\frac{(-1)^n}{n!}(-n)_k(\a+k+1)_{n-k},\;k=0,1,2,\ldots,n\n\eea
for all $n\in\{0,1,2,\ldots\}$.

By taking $y=x$ in (\ref{algJac}) we easily obtain
\bea\la{inv}& &\sum_{k=0}^n
\frac{(\a+\b+2k+1)(\a+\b+1)_k}{\G(\a+\b+n+k+2)}\times{}\nn
& &{}\hspace{1cm}{}\times P_k^{(\a,\b)}(x)P_{n-k}^{(-n-\a-1,-n-\b-1)}(x)=
\l\{\ba{ll}\ds\frac{1}{\G(\a+\b+1)}, & n=0\\ \\0, & n=1,2,3,\ldots\ea\r.\eea
for all $\a$ and $\b$. If we take $n=i-j$ in (\ref{inv}) and shift the
summation index we find
\bea & &\sum_{k=j}^i
\frac{(\a+\b+2k-2j+1)(\a+\b+1)_{k-j}}{\G(\a+\b+i-2j+k+2)}\times{}\nn
& &{}\hspace{1cm}{}\times
P_{i-k}^{(-i+j-\a-1,-i+j-\b-1)}(x)P_{k-j}^{(\a,\b)}(x)=
\frac{\delta_{ij}}{\G(\a+\b+1)},\;j\le i,\;i,j=0,1,2,\ldots.\n\eea
For $\a$ and $\b$ real with $\a+\b+1>-1$ we now obtain (\ref{invJac}) by
shifting both $\a$ and $\b$ by $j$.

Note that (\ref{algJac}) for $y=-x$ in a similar way leads to
\bea & &\sum_{k=j}^i
\frac{\a+\b+2k+1}{(\a+\b+k+j+1)_{i-j+1}}\times{}\nn
& &{}\hspace{1cm}{}\times
P_{i-k}^{(-\a-i-1,-\b-i-1)}(-x)P_{k-j}^{(\a+j,\b+j)}(x)=
\frac{x^{i-j}}{(i-j)!},\;j\le i,\;i,j=0,1,2,\ldots.\n\eea
This formula was used in \cite{Madrid}.

In order to prove (\ref{algJac}) we start with the left-hand side, apply
definition (\ref{defJac1}) to $P_k^{(\a,\b)}(x)$ and definition
(\ref{defJac2}) to $P_{n-k}^{(-n-\a-1,-n-\b-1)}(y)$ and change the order of
summation to obtain
\bea & &\sum_{k=0}^n
\frac{(\a+\b+2k+1)(\a+\b+1)_k}{\G(\a+\b+n+k+2)}P_k^{(\a,\b)}(x)
P_{n-k}^{(-n-\a-1,-n-\b-1)}(y)\nn
&=&\sum_{k=0}^n\sum_{i=0}^k\sum_{j=0}^{n-k}(-1)^{n-k}
\frac{(\a+\b+2k+1)(\a+\b+1)_k}{\G(\a+\b+n+k+2)}\frac{(\a+\b+k+1)_i}{i!}
\frac{(\a+i+1)_{k-i}}{(k-i)!}\times{}\nn
& &{}\hspace{1cm}{}\times
\frac{(\a+\b+n+k-j+2)_j}{j!}\frac{(\a+k+1)_{n-k-j}}{(n-k-j)!}
\l(\frac{x-1}{2}\r)^i\l(\frac{y-1}{2}\r)^j\nn
&=&\sum_{i=0}^n\sum_{k=i}^n\sum_{j=0}^{n-k}(-1)^{n-k}\times{}\nn
& &{}\hspace{1cm}{}\times
\frac{(\a+\b+2k+1)(\a+\b+1)_{i+k}(\a+i+1)_{n-i-j}}
{\G(\a+\b+n+k-j+2)\,i!\,(k-i)!\,j!\,(n-k-j)!}
\l(\frac{x-1}{2}\r)^i\l(\frac{y-1}{2}\r)^j\nn
&=&\sum_{i=0}^n\sum_{k=0}^{n-i}\sum_{j=0}^{n-i-k}(-1)^{n-i-k}\times{}\nn
& &{}\hspace{1cm}{}\times
\frac{(\a+\b+2i+2k+1)(\a+\b+1)_{2i+k}(\a+i+1)_{n-i-j}}
{\G(\a+\b+n+i+k-j+2)\,i!\,k!\,j!\,(n-i-k-j)!}
\l(\frac{x-1}{2}\r)^i\l(\frac{y-1}{2}\r)^j\nn
&=&\sum_{i=0}^n\sum_{j=0}^{n-i}\sum_{k=0}^{n-i-j}(-1)^{n-i}
\l(\frac{x-1}{2}\r)^i\l(\frac{y-1}{2}\r)^j\times{}\nn
& &{}\hspace{1cm}{}\times
\frac{(\a+\b+2i+2k+1)(\a+\b+1)_{2i+k}(\a+i+1)_{n-i-j}(-n+i+j)_k}
{\G(\a+\b+n+i-j+k+2)\,i!\,j!\,k!\,(n-i-j)!}\nn
&=&\sum_{i=0}^n\sum_{j=0}^{n-i}(-1)^{n-i}
\frac{(\a+\b+1)_{2i}(\a+i+1)_{n-i-j}}{i!\,j!\,(n-i-j)!}
\l(\frac{x-1}{2}\r)^i\l(\frac{y-1}{2}\r)^j\times{}\nn
& &{}\hspace{1cm}{}\times\sum_{k=0}^{n-i-j}
\frac{(-n+i+j)_k(\a+\b+2i+1)_k}{\G(\a+\b+n+i-j+k+2)\,k!}(\a+\b+2i+2k+1),
\;\ndots.\n\eea
Now we will show that for all $b$ we have
\be\la{nulalg}\sum_{k=0}^n\frac{(-n)_k(b)_k}{\G(b+n+k+1)\,k!}(b+2k)=0,
\;n=1,2,3,\ldots.\ee
In order to prove this we use the well-known Vandermonde summation formula
$$\hyp{2}{1}{-n,b}{c}{1}=\frac{(c-b)_n}{(c)_n},\;
(c)_n\ne 0,\;\ndots,$$
which can be written in a more general form as
$$\sum_{k=0}^n\frac{(-n)_k(b)_k}{\G(c+k)\,k!}=\frac{(c-b)_n}{\G(c+n)},
\;\ndots.$$
This formula is valid for all $b$ and $c$. By using this we find that for
all $b$ we have
\bea & &\sum_{k=0}^n\frac{(-n)_k(b)_k}{\G(b+n+k+1)\,k!}(b+2k)\nn
&=&b\,\sum_{k=0}^n\frac{(-n)_k(b+1)_k}{\G(b+n+k+1)\,k!}
-nb\,\sum_{k=0}^{n-1}\frac{(-n+1)_k(b+1)_k}{\G(b+n+k+2)\,k!}\nn
&=&b\,\frac{(n)_n}{\G(b+2n+1)}-
nb\,\frac{(n+1)_{n-1}}{\G(b+2n+1)}=0,\;n=1,2,3,\ldots,\n\eea
which proves (\ref{nulalg}). Now we use (\ref{nulalg}) to obtain
\bea & &\sum_{k=0}^n
\frac{(\a+\b+2k+1)(\a+\b+1)_k}{\G(\a+\b+n+k+2)}P_k^{(\a,\b)}(x)
P_{n-k}^{(-n-\a-1,-n-\b-1)}(y)\nn
&=&\sum_{i=0}^n(-1)^{n-i}\frac{(\a+\b+1)_{2i}}{i!\,(n-i)!}
\l(\frac{x-1}{2}\r)^i\l(\frac{y-1}{2}\r)^{n-i}\frac{\a+\b+2i+1}{\G(\a+\b+2i+2)}\nn
&=&\frac{1}{\G(\a+\b+1)}\frac{1}{n!}\sum_{i=0}^n{n \ch i}
\l(\frac{x-1}{2}\r)^i\l(\frac{1-y}{2}\r)^{n-i}\nn
&=&\frac{1}{\G(\a+\b+1)}\frac{1}{n!}\l(\frac{x-y}{2}\r)^n,\;\ndots,\n\eea
which proves (\ref{algJac}).

\section{Some remarks}

Note that we have from definition (\ref{defLag}) for the Laguerre
polynomials that
\be\la{monLag}L_n^{(-n)}(x)=(-1)^n\frac{x^n}{n!},\;\ndots.\ee
Hence, the polynomial $L_n^{(-n)}(x)$ reduces to a monomial of degree $n$
for all $n\in\{0,1,2,\ldots\}$. Definition (\ref{defJac3}) for the Jacobi
polynomials leads to
$$P_n^{(-n,\b)}(x)={n+\b \ch n}\l(\frac{x-1}{2}\r)^n,\;\ndots,$$
which is also a monomial. However, this monomial might reduce to the zero
polynomial. For instance, $P_n^{(-n,-n)}(x)$ equals the zero polynomial for
all $n\in\{1,2,3,\ldots\}$.

It is possible to generalize the Laguerre inversion formula (\ref{invLag}) to
\be\la{geninvLag}\sum_{k=0}^nL_k^{(\a+p_n)}(x)L_{n-k}^{(-\a-q_n)}(-x)=
\frac{(p_n-q_n+2)_n}{n!},\;\ndots,\ee
where $p_n$ and $q_n$ are arbitrary and even may depend on $n$. In order to
have an inversion formula we have to choose $p_n$ and $q_n$ such that
$$(p_n-q_n+2)_n=0,\;n=1,2,3,\ldots,$$
hence
$$p_n-q_n\in\{-n-1,-n,\ldots,-3,-2\},\;n=1,2,3,\ldots.$$
Note that the endpoint-cases $p_n-q_n=-n-1$ and $p_n-q_n=-2$ correspond to
the earlier mentioned inversion formulas (\ref{invLag}) and (\ref{invLag*})
respectively.

To prove (\ref{geninvLag}) we use (\ref{genLag}) to obtain
$$\L=\frac{1}{n!}\l.D_t^n\l[(1-t)^{-\a-1}
\exp\l(\frac{xt}{t-1}\r)\r]\r|_{t=0},\;\ndots,$$
where $\ds D_t=\frac{d}{dt}$ denotes differentiation with respect to $t$.
Hence by using Leibniz' rule we find
\bea & &\sum_{k=0}^nL_k^{(\a+p_n)}(x)L_{n-k}^{(-\a-q_n)}(-x)\nn
&=&\sum_{k=0}^n\frac{1}{k!}\l.D_t^k\l[(1-t)^{-\a-p_n-1}
\exp\l(\frac{xt}{t-1}\r)\r]\r|_{t=0}\times{}\nn
& &{}\hspace{2cm}{}\times\frac{1}{(n-k)!}
\l.D_t^{n-k}\l[(1-t)^{\a+q_n-1}\exp\l(\frac{-xt}{t-1}\r)\r]\r|_{t=0}\nn
&=&\frac{1}{n!}\l.D_t^n\l[(1-t)^{q_n-p_n-2}\r]\r|_{t=0}=
\frac{(p_n-q_n+2)_n}{n!},\;\ndots,\n\eea
which proves (\ref{geninvLag}).

Further we remark that if we replace $x$ by $\ds 1-\frac{2x}{\b}$ and $y$ by
$\ds 1-\frac{2y}{\b}$ in (\ref{algJac}), multiply by $\G(\a+\b+1)\b^n$ and
let $\b$ tend to infinity in an appropriate way we obtain by using
(\ref{limit})
\be\la{algLag}\sum_{k=0}^nL_k^{(\a)}(x)L_{n-k}^{(-n-\a-1)}(-y)=
\frac{(y-x)^n}{n!},\;\ndots.\ee
Note that (\ref{formLag}) is a special case of (\ref{algLag}) since
$$L_{n-k}^{(-n-\a-1)}(0)=(-1)^{n-k}\frac{(\a+k+1)_{n-k}}{(n-k)!}=
(-1)^n(-1)^k{n+\a \ch n-k},\;k=0,1,2,\ldots,n$$
for all $n\in\{0,1,2,\ldots\}$. Moreover, note that (\ref{algLag}) is a
special case of the well-known convolution formula for the classical Laguerre
polynomials
$$\sum_{k=0}^nL_k^{(\a)}(x)L_{n-k}^{(\b)}(y)=L_n^{(\a+\b+1)}(x+y),\;\ndots$$
in view of (\ref{monLag}). By using the technique demonstrated above this
convolution formula can be proved for all $\a$ and $\b$ which might even
depend on $n$.

Finally we remark that, by using the fact that
$$(b/2)_k(b+2k)=b(b/2+1)_k,\;k=0,1,2,\ldots,$$
formula (\ref{nulalg}) can also be obtained by using a summation formula for
a terminating well-poised hypergeometric series (see for instance formula
(III.9) in \cite{Slater}).

\vspace{10mm}

Menelaoslaan 4, 5631 LN Eindhoven, The Netherlands

\vspace{5mm}

Delft University of Technology, Faculty of Technical Mathematics and
Informatics,

P.O. Box 5031, 2600 GA Delft, The Netherlands, e-mail :
koekoek@twi.tudelft.nl

\end{document}